\documentclass[11pt]{article}
\usepackage{amsmath,amssymb,geometry,enumerate,float,cite,setspace}
\usepackage{color}
\newtheorem{theorem}{Theorem}
\newtheorem{lemma}[theorem]{Lemma}
\newtheorem{corollary}[theorem]{Corollary}

\begin{document}

\title{Induced $2$-Regular Subgraphs in $k$-Chordal Cubic Graphs}

\author{Michael A. Henning$^1$, Felix Joos$^2$, Christian L\"{o}wenstein$^2$, and Dieter Rautenbach$^2$}

\date{}

\maketitle

\begin{center}
{\small
$^1$ Department of Mathematics, University of Johannesburg, Auckland Park 2006, \\
South Africa, mahenning@uj.ac.za\\[3mm]
$^2$ Institute for Optimization and Operations Research, Ulm University, Ulm, \\
Germany, $\{$felix.joos,christian.loewenstein,dieter.rautenbach$\}$@uni-ulm.de
}
\end{center}

\begin{abstract}
We show that a cubic graph $G$ of order $n$ 
has an induced $2$-regular subgraph of order at least
\begin{itemize}  
\item $\frac{n-2}{4-\frac{4}{k}}$,
if $G$ has no induced cycle of length more than $k$,
\item $\frac{5n+6}{8}$,
if $G$ has no induced cycle of length more than $4$, and $n>6$, and
\item $\left(\frac{1}{4}+\epsilon\right)n$,
if the independence number of $G$
is at most $\left(\frac{3}{8}-\epsilon\right)n$.
\end{itemize}  
To show the second result
we give a precise structural description of cubic $4$-chordal graphs.
\end{abstract}



\section{Introduction}

The problem of finding a largest induced regular subgraph of a given graph goes back to Erd\H{o}s, Fajtlowicz, and Staton \cite{er}. It follows immediately from Ramsey's theorem \cite{ramsey} that every graph $G$ of order $n(G)$
has an induced regular subgraph of order $\Omega(\log n(G))$.
Special cases with fixed regularity such
as the independent set problem
or the induced matching problem have received a lot of attention.
In general, it is NP-hard to find a maximum induced (bipartite) $k$-regular subgraph of a given graph
as shown by Cardoso et al. \cite{cakalo},
who also extend the Hoffman upper bound on the independence number
to the maximum order of an induced $k$-regular subgraph.
Efficient algorithms for special graph classes \cite{lomopu},
exact exponential time algorithms \cite{gurasa},
as well as fpt-algorithms \cite{moth} for this problem have been studied.

\begin{figure}[H]
\begin{center}
\unitlength 1mm 
\linethickness{0.4pt}
\ifx\plotpoint\undefined\newsavebox{\plotpoint}\fi 
\begin{picture}(51,25)(0,0)
\put(0,5){\circle*{1.5}}
\put(0,15){\circle*{1.5}}
\put(10,15){\circle*{1.5}}
\put(20,15){\circle*{1.5}}
\put(30,15){\circle*{1.5}}
\put(40,15){\circle*{1.5}}
\put(50,15){\circle*{1.5}}
\put(10,5){\circle*{1.5}}
\put(20,5){\circle*{1.5}}
\put(30,5){\circle*{1.5}}
\put(40,5){\circle*{1.5}}
\put(50,5){\circle*{1.5}}
\put(0,15){\line(1,0){10}}
\put(10,15){\line(1,0){10}}
\put(20,15){\line(1,0){10}}
\put(30,15){\line(1,0){10}}
\put(40,15){\line(1,0){10}}
\put(10,5){\line(-1,0){10}}
\put(20,5){\line(-1,0){10}}
\put(30,5){\line(-1,0){10}}
\put(40,5){\line(-1,0){10}}
\put(50,5){\line(-1,0){10}}
\put(0,15){\line(1,-1){10}}
\put(10,15){\line(-1,-1){10}}
\put(20,15){\line(1,-1){10}}
\put(30,15){\line(-1,-1){10}}
\put(39,15){\line(1,0){1}}
\put(40,15){\line(1,-1){10}}
\put(50,15){\line(-1,-1){10}}
\qbezier(50,5)(25,-5)(0,5)
\qbezier(0,15)(25,25)(50,15)
\end{picture}
\end{center}
\caption{A graph $G$ with $c_{\rm ind}(G)=\frac{n(G)}{2}$.}\label{fign2}
\end{figure}
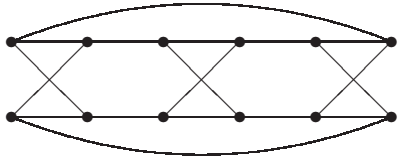

While the components induced by independent sets or induced matchings are clearly of bounded order,
there is no upper bound on the order of a component of an induced $k$-regular subgraph
for every $k$ at least $2$.
Unfortunately, this means that local techniques, 
which were successfully applied to independent sets and induced matchings,
hardly generalize to values of $k$ at least $2$.
Recently, Henning et al. \cite{hejolosa} studied 
the maximum order $c_{\rm ind}(G)$ of an induced $2$-regular subgraph of a given graph $G$.
They establish NP-hardness of $c_{\rm ind}(G)$ for graphs of maximum degree $4$.
For an $r$-regular graph $G$, they show
$$c_{\rm ind}(G)\geq \frac{n(G)}{2(r-1)}+\frac{1}{(r-1)(r-2)},$$
which implies $c_{\rm ind}(G)\geq \frac{n(G)+2}{4}$ if $G$ is cubic.
For a claw-free cubic graph $G$, they prove the asymptotically best-possible bound $c_{\rm ind}(G)>13n(G)/20$.
Furthermore, they believe that their general bound can be improved.
Specifically, for a cubic graph $G$, they conjecture $c_{\rm ind}(G)\geq \frac{n(G)}{2}$,
which would be best-possible in view of the graph in Figure \ref{fign2}.

In the present paper we study $c_{\rm ind}(G)$ for cubic graphs that do not have long induced cycles
or whose independence number is small.

For an integer $k$ at least $3$,
a graph $G$ is $k$-chordal if it does not have an induced cycle of length more than $k$.
Chordal graphs coincide with $3$-chordal graphs,
and graphs of small chordality were studied in \cite{both,chlosu,chrusr,krmanasa}.
Note that the components of an induced $2$-regular subgraph
of a $k$-chordal graph are of order at most $k$; that is, imposing $k$-chordality as an additional hypothesis
allows us to apply more local arguments.
Our results are an improvement of the bound from \cite{hejolosa}
for cubic $k$-chordal graphs
as well as a best-possible bound for $4$-chordal cubic graphs.
In order to prove this last result,
we give a precise structural description of $4$-chordal cubic graphs.

Before we proceed to our results and their proofs,
we would like to mention some further related notions and conjectures.

A set $D$ of vertices of a graph $G$ is a fair dominating set
if every vertex in $V(G)\setminus D$
has the same positive number of neighbors in $D$ \cite{cahahe}.
This definition implies that, if $G$ is an $r$-regular graph,
then a set $D$ of vertices of $G$ is a fair dominating set of $G$
if and only if $G-D$ is $s$-regular for some $s<r$.
Caro et al. \cite{cahahe} studied bounds on the fair domination number,
which is the minimum cardinality of a fair dominating set.
Clearly, for an $r$-regular graph $G$,
the fair domination number of $G$ is equal to
$n(G)-\max\{ n(H):\mbox{$H$ is an induced $s$-regular subgraph of $G$ with $s<r$}\}.$

Instead of regular induced subgraphs, one might consider induced subgraphs whose components are regular
but are allowed to have different degrees.
We conjecture that every cubic graph $G$
has an induced subgraph $H$ of order at least $\frac{3}{5}n(G)$
that is the disjoint union of $K_1$s, $K_2$s, and induced cycles.
The Petersen graphs shows that this is best possible.
\section{Results}

For a graph $G$, let $\kappa(G)$ be the number of components of $G$.
Recall that the cyclomatic number $\mu(G)$ of $G$ is $m(G)+\kappa(G)-n(G)$,
and that $G$ has a cycle if and only if $\mu(G)>0$. 
For a set $S$ of vertices of $G$, 
the closed neighborhood $N_G[S]$ of $S$ in $G$ 
contains $S$ and all neighbors of vertices in $S$.

\begin{lemma}\label{lemma1}
If $G$ is a connected cubic graph, then $G$ has an induced $2$-regular subgraph
with components $C_1,\ldots,C_t$ such that 
$$
\mu(G-V_{\le i-1})-\mu(G-V_{\le i})
\leq
\left\{
\begin{array}{ll}
2n(C_1) & \mbox {if $i=1$, and }\\
2n(C_i)-2 & \mbox {if $2\leq i\leq t$}.
\end{array}
\right.
$$
where $V_{\le 0} = \emptyset$ and $V_{\le i} = N_G[V(C_1) \cup \cdots \cup V(C_i)]$ for each $i$ 
with $1 \le i \le t$.
\end{lemma}
{\it Proof:}
We construct a sequence $G_0,G_1,\ldots,G_t$ of induced subgraphs of $G$
as well as a sequence $C_1,\ldots,C_t$ such that, 
for $i \in \{1,\ldots,t\}$, 
$C_i$ is an induced cycle of $G_{i-1}$,
and $G_i$ arises from $G_{i-1}$ by removing $C_i$ together with its neighbors;
that is, $G_{i}=G_{i-1}-V_i$ where $V_i=N_{G_{i-1}}[V(C_i)]$.
Clearly, $c_{\rm ind}(G)\geq \ell_1+\cdots+\ell_t$ where $\ell_i$ is the order of $C_i$; that is, $\ell_i=n(C_i)$.
Let $n_i=|V_i|$
and let $m_i$ be the number of edges of $G_{i-1}$ that are incident with a vertex in $V_i$;
that is,
$n_i=n(G_{i-1})-n(G_i)$ and $m_i=m(G_{i-1})-m(G_i)$.
Since $G$ is cubic, we have $n_i\leq 2\ell_i$.
Let $\mu_i=\mu(G_{i-1})-\mu(G_i)$.

Let $G_0=G$.
Let $C_1$ be any induced cycle of $G$.
If, for some $i\geq 2$, the graph $G_{i-1}$ has a cycle,
then choose $C_i$ as an induced cycle of $G_{i-1}$
such that $\mu_i-\ell_i$ is smallest possible.
The sequences terminate as soon as $G_i$ is a forest.
It remains to show that $\mu_1\leq 2\ell_1$ and that $\mu_i\leq 2\ell_i-2$ for $2\leq i\leq t$.

Let $i\in \{ 1,\ldots,t\}$.

If $V_i$ is the vertex set of a component of $G_{i-1}$, then
$\kappa(G_i)=\kappa(G_{i-1})-1$, and, since $G$ is cubic,
\begin{eqnarray*}
\mu_i
& = & (m(G_{i-1})-m(G_i))+(\kappa(G_{i-1})-\kappa(G_i))-(n(G_{i-1})-n(G_i))\\
& = & m_i+1-n_i\\
&\leq &\frac{3}{2}n_i+1-n_i\\
& = & \frac{1}{2}n_i+1\\
& \leq & \ell_i+1\\
& \leq & 2\ell_i-2.
\end{eqnarray*}
Hence, we may assume that
$V_i$ is not the vertex set of a component of $G_{i-1}$,
which implies that $\kappa(G_{i-1})-\kappa(G_i)\leq 0$.

Since $G$ is cubic, we have $m_i\leq \ell_i+3(n_i-\ell_i)=3n_i-2\ell_i$.
This implies
\begin{eqnarray*}
\mu_i
& \le & m_i - n_i\\
& \leq & 3n_i-2\ell_i-n_i\\
& = & 2n_i-2\ell_i\\
& \leq & 4\ell_i-2\ell_i\\
& = & 2\ell_i,
\end{eqnarray*}
which implies the desired bound for $i=1$.
Hence we may assume that $i\geq 2$.

If $C_i$ contains a vertex of degree $2$, then $n_i\leq 2\ell_i-1$,
and hence
\begin{eqnarray*}
\mu_i
& \leq & 2n_i-2\ell_i\\
& \leq & 4\ell_i-2-2\ell_i\\
& = & 2\ell_i-2.
\end{eqnarray*}
If $m_i$ is at most $3n_i-2\ell_i-2$, then a similar argument implies
$\mu_i\leq 2\ell_i-2$.

In view of the choice of $C_i$,
we may therefore assume that, for every induced cycle $C$ of $G_{i-1}$,
we have that $|V_C|=2n(C)$ where $V_C=N_ {G_{i-1}}[V(C)]$,
and that there are at least $4n(C)-1$ edges of $G_{i-1}$
that are incident with a vertex in $V_C$.
This implies that $V_C$ contains only vertices that are of degree~$3$ in $G_{i-1}$,
and that every vertex $v$ in $V_C\setminus C$
has at least one neighbor in $V(G_{i-1})\setminus V_C$.
Since $G_{i-1}$ is not a forest, it has a block $B$ that is distinct from $K_2$.
Since every vertex in $B$ lies on an induced cycle in $B$,
all vertices in $B$ have degree $3$ in $G_{i-1}$.
Since $i\geq 2$ and the graph $G$ is connected,
this implies that $B$ is not a component of $G_{i-1}$.
Let $u$ be a cutvertex of $G_{i-1}$.
Let $C_i$ be an induced cycle in $B$ that contains $u$.
Let $u^-$ and $u^+$ be the neighbors of $u$ in $C_i$.
Let $v$, $v^-$, and $v^+$ be the neighbors of $u$, $u^-$, and $u^+$ outside of $V(C_i)$, respectively.
Let $G'=G_{i-1}-(V_i\setminus \{ v,v^-,v^+\})$.
By construction, $v$ does not lie in the same component of $G'$ as $v^-$ or $v^+$.
If $v^-$ and $v^+$ lie in the same component of $G'$,
and $P$ is a shortest $v^-$-$v^+$-path in $G'$,
then $C'=(P\cup C_i)-\{ u\}$ is an induced cycle of $G_{i-1}$
with $|N_{G_{i-1}}[V(C')]|<2n(C')$, which is a contradiction.
Hence $v$, $v^-$, and $v^+$ all lie in different components of $G'$.
Since each of these vertices has a neighbor in $V(G_{i-1})\setminus V_i$,
we obtain $\kappa(G_{i-1})-\kappa(G_i)\leq -2$, and hence
\begin{eqnarray*}
\mu_i
& = & m_i+(\kappa(G_{i-1})-\kappa(G_i))-n_i\\
& \leq & 4\ell_i-2-2\ell_i\\
& = & 2\ell_i-2,
\end{eqnarray*}
which completes the proof. $\Box$

\begin{theorem}\label{theorem1}
If $G$ is a connected cubic $k$-chordal graph, then
$$c_{\rm ind}(G)\geq \frac{n(G)-2}{4-\frac{4}{k}}.$$
\end{theorem}
{\it Proof:} Let $C_1,\ldots,C_t$ be as in Lemma \ref{lemma1}.
We use the notation from the proof of Lemma \ref{lemma1}. By Lemma \ref{lemma1}, 
we have $\mu_1 \le 2 \ell_1$, and 
$\mu_i \le 2 \ell_i- 2$ for $2 \le i \le t$. 
Since $G$ is $k$-chordal, we have $\ell_i\leq k$ for all $i \ge 1$.
Therefore,
$$\mu_1\leq 2\ell_1=\frac{2\ell_1}{k}+\left(2-\frac{2}{k}\right)\ell_1\leq 2+\left(2-\frac{2}{k}\right)\ell_1$$
and, for $2 \le i \le t$, 
$$
\mu_i \leq  \left(2-\frac{2}{k}\right)\ell_i. 
$$
Since $G$ is a connected cubic graph, we have  $\mu(G)=\frac{n(G)}{2}+1$. Since $G_t$ is a forest, we have $\mu(G_t)=0$. Now
$$
\frac{n(G)}{2}+1
= \mu(G) 
= \left( \sum_{i=1}^t \mu_i \right) + \mu(G_t) =\sum_{i=1}^t \mu_i 
\le 2+\sum_{i=1}^t\left(2-\frac{2}{k}\right)\ell_i.
$$

This implies
$$c_{\rm ind}(G)\geq \sum\limits_{i=1}^t\ell_i\geq \frac{n(G)-2}{2\left(2-\frac{2}{k}\right)},$$
which completes the proof. $\Box$

\bigskip

\noindent It is obvious that the technique used in the proof of Lemma \ref{lemma1} and Theorem \ref{theorem1}
can also be applied to $r$-regular graphs for $r>3$.
Before we proceed to our result on $4$-chordal graphs,
we show another application of Lemma \ref{lemma1},
which relates $c_{\rm ind}(G)$ to the independence number $\alpha(G)$ of $G$.

\begin{theorem}\label{theorem1b}
Let $G$ be a connected cubic graph.
If $\alpha(G)\leq \left(\frac{3}{8}-\epsilon\right)n(G)$ for some $\epsilon>0$, then
$$c_{\rm ind}(G)>\left(\frac{1}{4}+\epsilon\right)n(G)-1.$$
\end{theorem}
{\it Proof:} For a contradiction, we suppose that $c_{\rm ind}(G)\leq \left(\frac{1}{4}+\epsilon\right)n(G)-1$.
Let $C_1,\ldots,C_t$ be as in Lemma \ref{lemma1}.
We use the notation from the proof of Lemma \ref{lemma1}.
Since $G_t$ is a forest,
$n_i\leq 2\ell_i$,
and no vertex of $C_1\cup \cdots \cup C_t$ is adjacent to a vertex of $G_t$,
we obtain
\begin{eqnarray*}
\alpha(G) & \geq & 
\sum\limits_{i=1}^t \left(  \frac{\ell_i-1}{2} \right) + \frac{n(G_t)}{2} \\
& = & \frac{1}{2} \left( 
\sum\limits_{i=1}^t \ell_i  + 
n(G) - \sum\limits_{i=1}^t n_i  \right)
- \frac{t}{2}  \\
& \ge & \frac{1}{2} \left( n(G) - 
\sum\limits_{i=1}^t \ell_i \right)
- \frac{t}{2}  \\
& \ge & \frac{1}{2} \left( n(G) - c_{\rm ind}(G) \right)
- \frac{t}{2}  \\
& > & \frac{1}{2}\left(n(G)-\left(\frac{1}{4}+\epsilon\right)n(G)\right)-\frac{t}{2}\\
& \geq & \left(\frac{3}{8}-\frac{\epsilon}{2}\right)n(G)-\frac{t}{2}.
\end{eqnarray*}
Together with $\alpha(G)\leq \left(\frac{3}{8}-\epsilon\right)n(G)$,
this implies $t\geq \epsilon n(G)$.
As in the proof of Theorem \ref{theorem1},
we obtain 
$$
\frac{n(G)}{2}+1=\sum_{i=1}^t\mu_i
\leq  2+\sum_{i=1}^t\left(2\ell_i-2\right)
= 2-2t+2\sum_{i=1}^t\ell_i
\leq 2-2\epsilon n(G)+2\sum_{i=1}^t\ell_i,
$$
which implies the contradiction
$c_{\rm ind}(G)\geq \sum\limits_{i=1}^t\ell_i\geq \left(\frac{1}{4}+\epsilon\right)n(G)-\frac{1}{2}$.~$\Box$

\bigskip

\noindent In order to prove our bound for $4$-chordal cubic graphs,
we describe their structure in detail.
Our next result characterizes all non-trivial blocks of a $4$-chordal cubic graph.
Let $K_n$, $P_n$, and $C_n$, be the complete graph, the path, and the cycle of order $n$, respectively.
Let $K_{n,m}$ be the complete bipartite graph with partite sets of order $n$ and $m$, respectively.
Let $G {\rm \, \Box \,} H$ be the Cartesian product of the graphs $G$ and $H$.

\begin{figure}[H]
\begin{center}
$\mbox{}$\hfill
\unitlength 1mm 
\linethickness{0.4pt}
\ifx\plotpoint\undefined\newsavebox{\plotpoint}\fi 
\begin{picture}(11,11)(0,0)
\put(5,0){\circle*{1.5}}
\put(10,5){\circle*{1.5}}
\put(0,5){\circle*{1.5}}
\put(5,10){\circle*{1.5}}
\put(0,5){\line(1,1){5}}
\put(5,10){\line(1,-1){5}}
\put(10,5){\line(-1,-1){5}}
\put(5,0){\line(-1,1){5}}
\put(0,5){\line(1,0){10}}
\end{picture}\hfill
\unitlength 1mm 
\linethickness{0.4pt}
\ifx\plotpoint\undefined\newsavebox{\plotpoint}\fi 
\begin{picture}(21,11)(0,0)
\put(5,0){\circle*{1.5}}
\put(10,5){\circle*{1.5}}
\put(0,5){\circle*{1.5}}
\put(5,10){\circle*{1.5}}
\put(0,5){\line(1,1){5}}
\put(5,10){\line(1,-1){5}}
\put(10,5){\line(-1,-1){5}}
\put(5,0){\line(-1,1){5}}
\put(0,5){\line(1,0){10}}
\put(20,5){\circle*{1.5}}
\put(5,10){\line(3,-1){15}}
\put(20,5){\line(-3,-1){15}}
\end{picture}\hfill
\unitlength 1mm 
\linethickness{0.4pt}
\ifx\plotpoint\undefined\newsavebox{\plotpoint}\fi 
\begin{picture}(14,11)(0,0)
\put(3,0){\circle*{1.5}}
\put(8,0){\circle*{1.5}}
\put(3,5){\circle*{1.5}}
\put(8,5){\circle*{1.5}}
\put(3,10){\circle*{1.5}}
\put(8,10){\circle*{1.5}}
\put(3,10){\line(0,-1){10}}
\put(8,10){\line(0,-1){10}}
\qbezier(3,0)(-3,5)(3,10)
\qbezier(8,10)(14,5)(8,0)
\put(8,10){\line(-1,0){5}}
\put(8,5){\line(-1,0){5}}
\put(8,0){\line(-1,0){5}}
\end{picture}\hfill
\unitlength 1mm 
\linethickness{0.4pt}
\ifx\plotpoint\undefined\newsavebox{\plotpoint}\fi 
\begin{picture}(11,11)(0,0)
\put(0,5){\circle*{1.5}}
\put(5,5){\circle*{1.5}}
\put(10,5){\circle*{1.5}}
\put(5,0){\circle*{1.5}}
\put(5,10){\circle*{1.5}}
\put(0,5){\line(1,1){5}}
\put(5,10){\line(1,-1){5}}
\put(10,5){\line(-1,-1){5}}
\put(5,0){\line(0,1){10}}
\put(0,5){\line(1,-1){5}}
\put(10,10){\circle*{1.5}}
\put(10,5){\line(0,1){5}}
\put(10,10){\line(-1,-1){5}}
\end{picture}\hfill$\mbox{}$
\end{center}
\caption{The graphs $D$, $D'$, $P_2 {\rm \, \Box \,} K_3$, and $K_{3,3}^-$.}\label{fig1}
\end{figure}
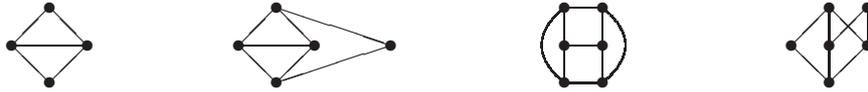

\begin{figure}[H]
\begin{center}
$\mbox{}$\hfill
\unitlength 1mm 
\linethickness{0.4pt}
\ifx\plotpoint\undefined\newsavebox{\plotpoint}\fi 
\begin{picture}(21,11)(0,0)
\put(0,0){\circle*{1.5}}
\put(0,10){\circle*{1.5}}
\put(0,10){\line(0,-1){10}}
\put(10,10){\line(1,0){2}}
\put(10,0){\line(1,0){2}}
\put(15,5){\makebox(0,0)[cc]{$\ldots$}}
\put(20,10){\circle*{1.5}}
\put(20,0){\circle*{1.5}}
\put(20,10){\line(0,-1){10}}
\put(10,10){\circle*{1.5}}
\put(10,0){\circle*{1.5}}
\put(0,10){\line(1,0){10}}
\put(10,10){\line(0,-1){10}}
\put(10,0){\line(-1,0){10}}
\put(18,10){\line(1,0){2}}
\put(18,0){\line(1,0){2}}
\end{picture}\hfill
\unitlength 1mm 
\linethickness{0.4pt}
\ifx\plotpoint\undefined\newsavebox{\plotpoint}\fi 
\begin{picture}(26,11)(0,0)
\put(0,5){\circle*{1.5}}
\put(5,0){\circle*{1.5}}
\put(5,10){\circle*{1.5}}
\put(5,10){\line(0,-1){10}}
\put(5,0){\line(-1,1){5}}
\put(0,5){\line(1,1){5}}
\put(15,10){\line(1,0){2}}
\put(15,0){\line(1,0){2}}
\put(20,5){\makebox(0,0)[cc]{$\ldots$}}
\put(25,10){\circle*{1.5}}
\put(25,0){\circle*{1.5}}
\put(25,10){\line(0,-1){10}}
\put(15,10){\circle*{1.5}}
\put(15,0){\circle*{1.5}}
\put(5,10){\line(1,0){10}}
\put(15,10){\line(0,-1){10}}
\put(15,0){\line(-1,0){10}}
\put(23,10){\line(1,0){2}}
\put(23,0){\line(1,0){2}}
\end{picture}\hfill
\unitlength 1mm 
\linethickness{0.4pt}
\ifx\plotpoint\undefined\newsavebox{\plotpoint}\fi 
\begin{picture}(31,11)(0,0)
\put(0,5){\circle*{1.5}}
\put(5,0){\circle*{1.5}}
\put(5,10){\circle*{1.5}}
\put(5,10){\line(0,-1){10}}
\put(5,0){\line(-1,1){5}}
\put(0,5){\line(1,1){5}}
\put(15,10){\line(1,0){2}}
\put(15,0){\line(1,0){2}}
\put(20,5){\makebox(0,0)[cc]{$\ldots$}}
\put(25,10){\circle*{1.5}}
\put(25,0){\circle*{1.5}}
\put(25,10){\line(0,-1){10}}
\put(30,5){\circle*{1.5}}
\put(25,10){\line(1,-1){5}}
\put(30,5){\line(-1,-1){5}}
\put(15,10){\circle*{1.5}}
\put(15,0){\circle*{1.5}}
\put(5,10){\line(1,0){10}}
\put(15,10){\line(0,-1){10}}
\put(15,0){\line(-1,0){10}}
\put(23,10){\line(1,0){2}}
\put(23,0){\line(1,0){2}}
\end{picture}\hfill$\mbox{}$
\end{center}
\caption{$B_k=P_2 {\rm \, \Box \,} P_k$, $B_k'$, and $B_k''$ for $k\geq 2$.}\label{fig2}
\end{figure}
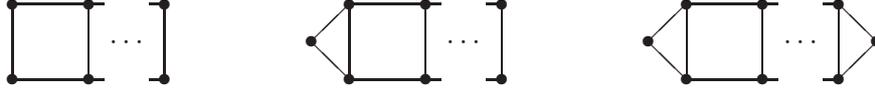

For some integer $k$ at least $2$, let $B_k$ be the graph $P_2 {\rm \, \Box \,} P_k$.
Note that $B_2$ is $C_4$.
Let $B_k'$ arise from $B_k$ by adding a new vertex to $B_k$ and joining it to two adjacent vertices of $B_k$ of degree $2$.
Let $B_k''$ arise from $B_k'$ by adding a new vertex to $B_k'$ and joining it to the two adjacent vertices of $B_k'$ of degree~$2$.
See Figure \ref{fig2} for an illustration.

Let
\begin{eqnarray*}
{\cal F}&=&\{ K_3,K_4,D,D',P_2 {\rm \, \Box \,} K_3,K_{2,3},K_{3,3},K_{3,3}^-\}\mbox{, and}\\
{\cal B}&=&\{ B_k:k\geq 2\}\cup \{ B_k':k\geq 2\}\cup \{ B_k'':k\geq 2\}.
\end{eqnarray*}

\begin{theorem} \label{theorem2}
If $G$ is a $2$-connected subcubic $4$-chordal graph, then $G$ belongs to ${\cal F}\cup {\cal B}$.
\end{theorem}
{\it Proof:} Let $G$ be a $2$-connected subcubic $4$-chordal graph.
Since all graphs in ${\cal F}$ are $2$-connected, subcubic, and $4$-chordal,
we may assume that $G$ does not belong to ${\cal F}$.
We consider different cases.

First, we assume that $G$ contains the diamond $D$ as an induced subgraph.
Let $a$ and $b$ be the two vertices of degree $2$ in $D$.
Since $G$ is not $D$, we may assume that $a$ has a neighbor $c$ not in $D$.
Since $G$ is not $D'$, the vertex $b$ is not adjacent to $c$.
Since $G$ is $2$-connected, there is a shortest path $P$ between $b$ and $c$
that does not intersect $D-b$.
Now $P$ together with a shortest $a$-$b$-path in $D$ yields an induced cycle of length more than $4$,
which is a contradiction.
Therefore, we may assume that $G$ is $D$-free.

Next, we assume that $G$ contains a triangle $T:abca$.
If all vertices of $T$ have degree $3$ in $G$,
then, since $G$ is $D$-free and not $K_4$,
the vertices $a$, $b$, and $c$ have distinct neighbors, say $a'$, $b'$, and $c'$, outside of $T$, respectively.
Since $G$ is $4$-chordal, no two of the vertices $a'$, $b'$, and $c'$
are joined by an induced path of length at least $2$ in $G-V(T)$.
Since $G$ is $2$-connected, every two of the vertices $a'$, $b'$, and $c'$
are joined by an induced path in $G-V(T)$.
Hence $a'$, $b'$, and $c'$ induce a $K_3$, and $G$ is $P_2 {\rm \, \Box \,} K_3$, which is a contradiction.
Therefore, we may assume that $b$ has degree $2$ in $G$.
Since $G$ is not $K_3$, we may assume that $a$ has degree $3$ in $G$.
Let $a_1$ be the neighbor of $a$ outside of $T$.
Since $G$ is $2$-connected, the graph $G-\{ a,b\}$ contains a shortest $a_1$-$c$-path $P$.
Since $G$ is $4$-chordal and $D$-free, the path $P$ has order exactly $3$.
Let $c_1$ be the unique internal vertex of $P$.
Note that $G[\{ a,b,c,a_1,c_1\}]$ is isomorphic to $B_2'$.
If, for some $k\geq 2$, a proper induced subgraph $G'$ of $G$ is isomorphic to $B_k'$,
and $a_k$ and $c_k$ are the two adjacent vertices of degree $2$ in $G'$,
then we may assume that $a_k$ has a neighbor $a_{k+1}$ outside of $G'$.
Since $G$ is $2$-connected, the graph $G-(V(G')\setminus \{ c_k\})$
contains a shortest $a_{k+1}$-$c_k$-path $Q$.
Since $G$ is $4$-chordal, the path $Q$ has order at most $3$.
If $Q$ has order $2$, then $G$ is $B_k''$.
If $Q$ has order $3$, then $G$ has an induced subgraph that is isomorphic to $B_{k+1}'$.
By an inductive argument, we obtain that $G$ is $B_k'$ or $B_k''$ for some $k\geq 2$.
Therefore, we may assume that $G$ is triangle-free.

Next, we assume that $G$ contains $K_{2,3}$ as an induced subgraph.
Let $a$, $b$, and $c$ be the three vertices of degree $2$ in this $K_{2,3}$.
Since $G$ is not $K_{2,3}$, we may assume that $a$ has a neighbor $a'$ outside of $K_{2,3}$.
Since $G$ is $2$-connected, we may assume, by symmetry between $b$ and $c$,
that $P$ is a shortest $a'$-$b$-path in $G-(V(K_{2,3})\setminus \{ b,c\})$.
Since $G$ is $4$-chordal, the path $P$ has order $2$;
that is, the vertex $a'$ is adjacent to $b$.
Since $G$ is not $K_{3,3}$, the vertex $a'$ is not adjacent to $c$.
Since $G$ is not $K_{3,3}^-$, we may assume that $Q$ is a shortest $a'$-$c$-path in $G-(V(K_{2,3})\setminus \{ c\})$.
Since $Q$ is of order at least $3$, it is contained in an induced cycle of length at least $5$ in $G$, which is a contradiction.
Therefore, we may assume that $G$ is $K_{2,3}$-free.

Since $G$ is $4$-chordal, $2$-connected, and triangle-free,
it contains an induced $4$-cycle $C:a_1a_2b_2b_1a_1$.
Since $C$ is isomorphic to $B_2$, we may assume that $G$ is not $C$.
Therefore, we may assume, by symmetry, that $a_2$ has a neighbor $a_3$ outside of $C$.
Since $G$ is $2$-connected, we may assume that $P$ is a shortest path in $G-a_2$
between $a_3$ and a vertex in $\{ a_1,b_1,b_2\}$.
If $P$ is an $a_3$-$b_1$-path, then, since $G$ is triangle-free and $4$-chordal,
the path $P$ has order $2$;
that is, the vertex $a_3$ is adjacent to $b_1$,
and $G$ is not $K_{2,3}$-free, which is a contradiction.
Hence, we may assume, by symmetry between $a_1$ and $b_2$,
that $P$ is an $a_3$-$b_2$-path.
Since $G$ is triangle-free and $4$-chordal, the path $P$ has order exactly $3$.
Let $b_3$ be the unique internal vertex of $P$.
Since $G$ is triangle-free and $K_{2,3}$-free,
the graph $G[\{ a_1,a_2,a_3,b_1,b_2,b_3\}]$ is isomorphic to $B_3$.
Applying an inductive argument as above, we obtain that $G$ is $B_k$ for some $k\geq 3$,
which completes the proof. $\Box$

\bigskip

\noindent Let $G$ be a cubic $4$-chordal graph.
If $B$ is a block of $G$ that is distinct from $K_2$,
then Theorem \ref{theorem2} 
implies that $B$ belongs to ${\cal F}\cup {\cal B}$.
Furthermore, all edges of $G$ between a vertex in $V(B)$ and a vertex in $V(G)\setminus V(B)$ are bridges.
Considering the vertex degrees of the graphs in ${\cal F}\cup {\cal B}$,
this implies that there are at most four edges between $V(B)$ and $V(G)\setminus V(B)$.
Therefore, contracting every block of $G$ that is distinct from $K_2$ to a single vertex,
results in a tree of maximum degree $4$.
Reversing this process leads to the following constructive description of cubic $4$-chordal graphs.

\begin{corollary}\label{corollary1}
If $G$ is a connected cubic $4$-chordal graph, then $G$ is either $2$-connected,
in which case $G$ belongs to $\{ K_4,K_{3,3},P_2 {\rm \, \Box \,} K_3\}$,
or $G$ arises from a tree $T$ of order at least $2$ and maximum degree at most $4$ by replacing
\begin{itemize}
\item every endvertex of $T$ with $D'$,
\item every vertex of $T$ of degree $2$ with either $D$, or $K_{3,3}^-$, or $B_k''$ for some $k\geq 2$,
\item some vertices of $T$ of degree $3$ with either $K_3$, or $K_{2,3}$, or $B_k'$ for some $k\geq 2$, and
\item every vertex of $T$ of degree $4$ with $B_k$ for some $k\geq 2$.
\end{itemize}
\end{corollary}
Based on this structural description,
we proceed to our second main result.

\begin{theorem}\label{theorem3}
If $G$ is a connected cubic $4$-chordal graph that does not belong to the set $\{ K_4,K_{3,3},P_2 {\rm \, \Box \,} K_3\}$,
then
$$c_{\rm ind}(G)\geq \frac{5}{8}n(G)+\frac{3}{4}.$$
Furthermore, equality holds if and only if
$G$ arises from a tree $T$ of order at least $2$ and maximum degree at most $3$ by replacing
\begin{itemize}
\item every endvertex of $T$ with $D'$,
\item every vertex of $T$ of degree $2$ with $B_3''$, and
\item every vertex of $T$ of degree $3$ with $K_3$.
\end{itemize}
\end{theorem}
{\it Proof:} We prove the statement by induction on the order of $G$.
By Corollary \ref{corollary1}, the graph $G$ is not $2$-connected.
Let $T$ be as in the statement of Corollary \ref{corollary1}.
Since $T$ has at least two endvertices, the order of $G$ is at least $10$.
If $n(G)=10$, then $G$ arises from the disjoint union of two copies of $D'$
by connecting the two vertices of degree $2$ by a bridge, and $c_{\rm ind}(G)=7=\frac{5}{8}n(G)+\frac{3}{4}$.
Now let $n(G)>10$, which implies that $T$ has order at least $3$.

\begin{figure}[H]
\begin{center}
$\mbox{}$\hfill
\unitlength 1mm 
\linethickness{0.4pt}
\ifx\plotpoint\undefined\newsavebox{\plotpoint}\fi 
\begin{picture}(41,19)(0,0)

\put(0,15){\line(-1,0){5}}
\put(0,5){\line(-1,0){5}}

\put(40,15){\line(1,0){5}}
\put(40,5){\line(1,0){5}}

\put(0,5){\circle*{1.5}}
\put(20,5){\circle*{1.5}}
\put(0,15){\circle*{1.5}}
\put(20,15){\circle*{1.5}}
\put(0,15){\line(0,-1){10}}
\put(20,15){\line(0,-1){10}}
\put(10,15){\circle*{1.5}}
\put(30,15){\circle*{1.5}}
\put(10,5){\circle*{1.5}}
\put(30,5){\circle*{1.5}}
\put(0,15){\line(1,0){10}}
\put(20,15){\line(1,0){10}}
\put(10,15){\line(0,-1){10}}
\put(30,15){\line(0,-1){10}}
\put(10,5){\line(-1,0){10}}
\put(30,5){\line(-1,0){10}}
\put(10,15){\line(1,0){10}}
\put(20,5){\line(-1,0){10}}
\put(40,15){\circle*{1.5}}
\put(40,5){\circle*{1.5}}
\put(30,15){\line(1,0){10}}
\put(40,15){\line(0,-1){10}}
\put(40,5){\line(-1,0){10}}
\put(0,19){\makebox(0,0)[cc]{$a_1$}}
\put(0,1){\makebox(0,0)[cc]{$b_1$}}
\put(10,19){\makebox(0,0)[cc]{$a_2$}}
\put(10,1){\makebox(0,0)[cc]{$b_2$}}
\put(20,19){\makebox(0,0)[cc]{$a_3$}}
\put(20,1){\makebox(0,0)[cc]{$b_3$}}
\put(30,19){\makebox(0,0)[cc]{$a_4$}}
\put(30,1){\makebox(0,0)[cc]{$b_4$}}
\put(40,19){\makebox(0,0)[cc]{$a_5$}}
\put(40,1){\makebox(0,0)[cc]{$b_5$}}
\end{picture}
\hfill
\unitlength 1mm 
\linethickness{0.4pt}
\ifx\plotpoint\undefined\newsavebox{\plotpoint}\fi 
\begin{picture}(41,19)(0,0)
\put(0,15){\line(-1,0){5}}
\put(0,5){\line(-1,0){5}}

\put(40,15){\line(1,0){5}}
\put(40,5){\line(1,0){5}}

\put(0,5){\circle*{1.5}}
\put(0,15){\circle*{1.5}}
\put(0,15){\line(0,-1){10}}
\put(0,15){\line(1,0){10}}
\put(10,5){\line(-1,0){10}}
\put(40,15){\circle*{1.5}}
\put(40,5){\circle*{1.5}}
\put(30,15){\line(1,0){10}}
\put(40,15){\line(0,-1){10}}
\put(40,5){\line(-1,0){10}}
\put(0,19){\makebox(0,0)[cc]{$a_1$}}
\put(0,1){\makebox(0,0)[cc]{$b_1$}}
\put(40,19){\makebox(0,0)[cc]{$a_5$}}
\put(40,1){\makebox(0,0)[cc]{$b_5$}}
\put(30,15){\line(-1,0){22}}
\put(30,5){\line(-1,0){22}}
\end{picture}
\hfill$\mbox{}$
\end{center}
\caption{An induced $B_5$ in $G$ and $G''$.}\label{fig3}
\end{figure}
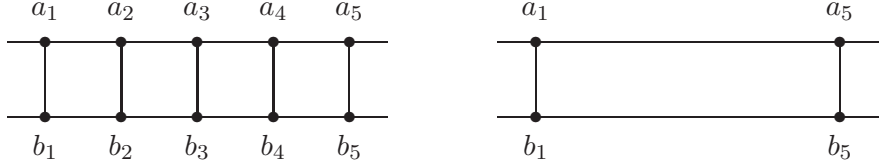

First, we assume that $G$ contains $B_5$ as an induced subgraph.
We denote its vertices as in the left of Figure \ref{fig3}.

Note that $G'=G-\{ a_2,a_3,a_4,b_2,b_3,b_4\}$ has exactly two components.
Let $G''$ arise from $G'$ by adding the two edges $a_1a_5$ and $b_1b_5$.
Clearly, $G''$ is cubic.
Since removing any two edges of $K_{3,3}$ or $P_2 {\rm \, \Box \,} K_3$ does not disconnect these graphs,
$G''$ does not belong to $\{ K_4,K_{3,3},P_2 {\rm \, \Box \,} K_3\}$.
Note that $a_1a_5b_5b_1a_1$ is the only induced cycle of $G''$ that is not also a cycle of $G$.
Therefore, $G''$ is $4$-chordal.
Let $H''$ be an induced $2$-regular subgraph of $G''$.
If $H''$ contains the cycle $a_1a_5b_5b_1a_1$,
then let $H=(H''-\{ a_1,a_5,b_5,b_1\})\cup a_1a_2b_2b_1a_1\cup a_4a_5b_5b_4a_4$.
If $H''$ does not contain the cycle $a_1a_5b_5b_1a_1$,
then $H''$ does not contain any of the two edges $a_1a_5$ and $b_1b_5$.
Therefore, if $a_1\in V(H'')$, then $a_5,b_5\not\in V(H'')$,
and let $H=H''\cup a_3a_4b_4b_3a_3$.
By symmetry, this implies in all cases
that $G$ has an induced $2$-regular subgraph $H$ with $n(H)\geq n(H'')+4$.
By induction, we obtain
\begin{eqnarray*}
c_{\rm ind}(G) & \geq & c_{\rm ind}(G'')+4
\geq \frac{5}{8}(n(G)-6)+\frac{3}{4}+4
 >  \frac{5}{8}n(G)+\frac{3}{4}.
\end{eqnarray*}
Therefore, we may assume that $G$ is $B_5$-free.
By a similar argument, we may assume that $G$ is $B_4'$-free.

Let $P:uvw\ldots$ be a longest path in $T$.
Note that $P$ has order at least $3$, and that all neighbors of $v$ in $T$ that are distinct from $w$ are endvertices of $T$.
Let $U=N_T(v)\setminus \{ w\}$.
The vertex $v$ in $T$ is either a vertex of $G$ that is of degree $3$
or it corresponds to a block $B$ in $G$ according to Corollary \ref{corollary1}.
In the first case, let $B=K_1$.
Every vertex in $U$ corresponds to an induced $D'$ in $G$ that is connected to $B$ by a bridge of $G$.
There is a unique vertex $x$ of $G$ that does not belong to $B$ or to one of the copies of $D'$ that correspond to the vertices in $U$, such that $x$ has a neighbor $y$ in $B$.
Note that $xy$ is a bridge of $G$.
Let $G'$ be the component of $G-xy$ that contains $x$,
and let $B^+$ be the component of $G-xy$ that contains $y$.
Let $G''$ arise from the disjoint union of $G'$ and $D'$ by adding an edge between $x$ and the vertex of degree $2$ in $D'$.
Note that $G''$ is a cubic $4$-chordal graph of order less than $G$ that does not belongs to $\{ K_4,K_{3,3},P_2 {\rm \, \Box \,} K_3\}$.
Let $H''$ be an induced $2$-regular subgraph of $G''$ of order $c_{\rm ind}(G'')$.
If $x\in V(H'')$, then $H''$ contains exactly three vertices of the $D'$ that was added to $G'$.
Therefore, $c_{\rm ind}(G)\geq c_{\rm ind}(G'')-3+c_{\rm ind}(B^+-y)$.
Similarly,
if $x\not\in V(H'')$, then $H''$ contains exactly four vertices of the $D'$ that was added to $G'$,
and hence $c_{\rm ind}(G)\geq c_{\rm ind}(G'')-4+c_{\rm ind}(B^+)$.
The following table summarizes relevant values for all possibilities for $B$.

$$
\begin{array}{|l|l|l|l|l|}\hline
B & c_{\rm ind}(B^+) & c_{\rm ind}(B^+-y) & c_{\rm ind}(G)-c_{\rm ind}(G'')\geq & \frac{5}{8}n(G)-\frac{5}{8}n(G'')\\ \hline\hline
D & 7 & 6 & 3 & 2.5\\\hline
K_{3,3}^-& 8 & 8 & 4 & 3.75\\\hline
B_2''& 8 & 8 & 4 & 3.75\\\hline
B_3''& 9 & 8 & 5 & 5\\\hline
K_1& 8 & 8 & 4 & 3.75\\\hline
K_3& 9 & 8 & 5 & 5\\\hline
K_{2,3} & 11 & 10 & 7 & 6.25\\\hline
B_2' & 11 & 10 & 7 & 6.25 \\\hline
B_3' & 12 & 12 & 8 & 7.5\\\hline
B_2 & 13 & 12 & 9 & 8.75\\\hline
B_3 & 15 & 14 & 11 & 10\\\hline
B_4 & 16 & 16 & 12 & 11.25 \\\hline
\end{array}
$$

\bigskip

\begin{figure}[H]
\begin{center}
\unitlength 1mm 
\linethickness{0.4pt}
\ifx\plotpoint\undefined\newsavebox{\plotpoint}\fi 
\begin{picture}(86,11)(0,0)
\put(30,5){\circle*{1.5}}
\put(70,5){\circle*{1.5}}
\put(35,0){\circle*{1.5}}
\put(75,0){\circle*{1.5}}
\put(35,10){\circle*{1.5}}
\put(75,10){\circle*{1.5}}
\put(75,0){\line(-1,1){5}}
\put(70,5){\line(1,1){5}}
\put(55,10){\circle*{1.5}}
\put(55,0){\circle*{1.5}}
\put(60,5){\circle*{1.5}}
\put(45,10){\circle*{1.5}}
\put(85,10){\circle*{1.5}}
\put(45,0){\circle*{1.5}}
\put(85,0){\circle*{1.5}}
\put(35,10){\line(1,0){10}}
\put(75,10){\line(1,0){10}}
\put(45,10){\line(0,-1){10}}
\put(85,10){\line(0,-1){10}}
\put(45,0){\line(-1,0){10}}
\put(85,0){\line(-1,0){10}}
\put(45,10){\line(1,0){10}}
\put(55,0){\line(-1,0){10}}
\put(75,10){\line(1,-1){10}}
\put(85,10){\line(-1,-1){10}}
\put(70,5){\line(-1,0){10}}
\put(12,5){\oval(24,10)[]}
\put(20,5){\circle*{1.5}}
\put(20,5){\line(1,0){10}}
\put(16,5){\makebox(0,0)[cc]{$x$}}
\put(30,9){\makebox(0,0)[cc]{$y$}}
\put(8,5){\makebox(0,0)[cc]{$G'$}}
\put(35,10){\line(-1,-1){5}}
\put(30,5){\line(1,-1){5}}
\put(35,0){\line(0,1){10}}
\put(55,10){\line(0,-1){10}}
\put(55,0){\line(1,1){5}}
\put(60,5){\line(-1,1){5}}
\end{picture}
\end{center}
\caption{The case $B=B_3''$.}\label{figb3ss}
\end{figure}
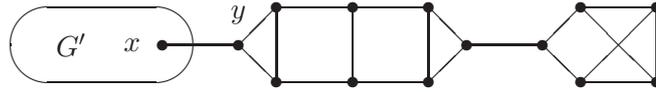

\noindent Figure \ref{figb3ss} illustrates the case $B=B_3''$. In this case, combining two triangles in $B$ with a triangle in $D'$ yields $c_{\rm ind}(B^+)=9$. Combining two cycles of length $4$, one in $B$ and one in $D'$, yields $c_{\rm ind}(B^+-y)=8$.

\noindent Note that for $B\in \{ B_2',B_3'\}$, there are two non-isomorphic configurations for $G$.
Since these lead to the same values, we do not distinguish them within the table.

Since the entries in the second to last column are consistently at least as large as the entries in the last column,
the desired bound follows by induction.
The statement about the extremal graphs easily follows from the base case of the induction,
and the fact that only $B_3''$ and $K_3$ lead to equal values within the last two columns of the table. $\Box$

\bigskip

\noindent {\bf Acknowledgements} 
The research of the first author was supported in part 
by the 
South African National Research Foundation and the 
University of Johannesburg.

\end{document}